\documentclass[10pt]{article}

\usepackage{amsfonts, latexsym, amsmath, amssymb, fancyhdr, amsthm}
\usepackage[english]{babel}

\newtheorem{theorem}{Theorem}[section]
\newtheorem{lemma}{Lemma}[section]
\newtheorem{proposition}{Proposition}[section]
\newtheorem{definition}{Definition}[section]
\newtheorem{corollary}{Corollary}[section]

\title{{\bf Global Stability of a Partial Differential Equation with Distributed Delay due to Cellular Replication\thanks{This paper has been published in Nonlinear Analysis, 54, 8, 1469-1491 (2003).}}}
\author{Mostafa Adimy\thanks{mostafa.adimy@univ-pau.fr} \quad and \quad Fabien Crauste\thanks{fabien.crauste@univ-pau.fr}}
\date{Year 2002}

\begin{document}

\maketitle

\begin{center}
{\em Laboratoire de Math\'ematiques Appliqu\'ees\\Universit\'e de
Pau et des Pays de l'Adour\\
Avenue de l'universit\'e, 64000 Pau, France}
\end{center}

\bigskip{}

\begin{abstract}
In this paper, we investigate a nonlinear partial differential equation, arising from a model of cellular proliferation. This model describes the production of blood cells in the bone marrow. It is represented by a partial differential equation with a retardation of the maturation variable and a distributed temporal delay. Our aim is to prove that the behaviour of primitive cells influences the global behaviour of the population.
\end{abstract}

\section{Introduction and motivation} \label{introduction}

This paper analyses a general model of the blood production system based on a model proposed by
Mackey and Rey in 1993 \cite{mackey1993}. The initial form of this model is a time-age-maturity
structured system and it describes the dynamics of proliferative stem cells and precursors in the
bone marrow. It consists in a population of cells which are capable of both proliferation and
maturation.

In this model, the period of life of each cell is divided into a resting phase and a proliferating
phase (see \cite{burns}). The cells in the resting phase can not divide. They mature and, provided
they do not die, they eventually enter the proliferating phase. In the proliferating phase, if they
do not die by apoptosis, the cells are committed to divide and give birth, at the point of
cytokinesis, to two daughter cells. The two daughter cells enter immediately the resting phase and
complete the cycle. In the resting phase, a cell can remains indefinitely.

The model in \cite{mackey1993} has been analysed by Mackey and Rey in 1995
\cite{mackey1995_1,mackey1995_2}, Crabb et al. in 1996 \cite{crabb1996_1,crabb1996_2}, Dyson,
Villella-bressan and Webb in 1996 \cite{webb1996} and Adimy and Pujo-Menjouet in 2001
\cite{adimypujo,adimypujo2}. In these studies, the model in \cite{mackey1993} was simplified by
assuming that all cells divide exactly at the same age. In the most general situation in a cellular
population, it is believed that the time required for a cell to divide is not identical between
cells (see the works of Bradford et al. \cite{bradford} about experiments on mice). To our
knowledge, this hypothesis has been given for the first time in \cite{mackey1993} for particular
cases, and only numerical studies have been investigated.

In \cite{webb2000} and \cite{webb2000_2}, Dyson et al. considered a time-age-maturity structured
equation in which the age for a cell to divide is not identical between cells. They presented the
basic theory of existence and uniqueness and properties of the solution operator. However, in their
model, the division is represented by the following boundary condition
\begin{equation} \label{bcwebb}
p(t,m,0)=\int_{0}^{+\infty} \beta(a)p(t,m,a)da,
\end{equation}
which considers only one phase (the proliferating one), and the intermediary flux between the two
phases is not represented in this model.

In \cite{mackey1994} (1994) and \cite{mackey1999} (1999), Mackey and Rudnicki studied the behaviour of solutions of the model considered in \cite{mackey1993}, but only in the case when all cells divide at the same age. They obtained a first order partial differential equation with discrete time delay $\tau$ and a nonlocal dependence in the maturity variable $\Delta(m)$ due to cell replication. That is,
\begin{equation}\label{eqmackey}
\frac{\partial u}{\partial t}(t,m)+V(m)\frac{\partial u}{\partial m}(t,m) = f(u(t,m),u(t-\tau,\Delta(m))),
\end{equation}
with $\Delta:[0,1]\to[0,1]$ a continuous function satisfying $\Delta(0)=0$ and $\Delta(m)<m$ for $m\in(0,1]$.

They gave in \cite{mackey1999}, which is more general than \cite{mackey1994}, a criterion for global stability in such equations. However, these authors considered only the special case when the term $f(u,v)$ in (\ref{eqmackey}) does not depend on the maturity variable. Moreover, the nonlocal function $\Delta(m)$ usually depends on $\tau$ and the condition $\Delta(m)<m$, used by Mackey and Rudnicki in the example $V(m)=rm$, $r>0$, and $g(m)=\alpha m$, $0<\alpha\leq 1$, which yields to $\Delta(m)=\alpha^{-1}e^{-r\tau}m$, is not true in the general case, for all $\tau>0$.

In this work, we consider a general situation, when the age at cytokinesis is distributed with a
kernel supported on an interval $[\underline{\tau},\overline{\tau}]$ with
$0\leq\underline{\tau}<\overline{\tau}<+\infty$. This yields to the boundary condition
(\ref{bcn(t,m,a)}).

We obtain a first order partial differential equation with a distributed time delay and a retardation of the maturation variable, which depends also on the time delay. In our model, the nonlinear part depends on the maturity variable (see Equations (\ref{N(t,m)}) and (\ref{P(t,m)})).

The population behaviour of the model obtained is analysed in the case
\begin{equation} \label{equationintro}
\lim_{t\to +\infty} N(t,0)=0.
\end{equation}
In the interpretation of this model, $N(t,0)$ corresponds to the population of the primitive cells
type, and Condition (\ref{equationintro}) corresponds, for a long time, to an insufficient supply
of the primitive cells.

We prove that Condition (\ref{equationintro}) guarantees, under local stability (resp. local
exponential stability) of the partial differential equation, global stability (resp. global
exponential stability) of the system. These results are obtained under the assumption that the
cells have enough time, in the proliferating phase ($\underline{\tau} > \tau_{0}$), to sufficiently
increase their maturity.

It is also proved that, if this minimum time of replication $\underline{\tau}$ is too large, then
the destruction of the population of stem cells (cells with small maturity) affects the entire
system, and the population extincts for a finite time (see Corollary \ref{extinction}).

The paper is organised as follows. In Section \ref{presentationmodel}, we present the model. We obtain a system of two partial differential equations with distributed time delay and a nonlocal dependence in the maturity variable. We give an integrated formulation of our model by using the classical variation of constant formula. In Section \ref{sectionlocalstability}, we study the local stability of the trivial solution. In Section \ref{stabilitydde}, we investigate the global stability of the delay differential equation satisfied by the population of the primitive cells type. We use a Liapunov functional to prove the global stability of the trivial solution of this equation. In Section \ref{globalstability}, we provide a criterion for global stability of the trivial solution in our model. This criterion is based on Condition  (\ref{equationintro}) and local stability of the trivial solution.

\section{Presentation of the model} \label{presentationmodel}

Each cell is characterised, in the two phases, by its age and its maturity. The maturity describes
the development of the cell. It is the concentration of what composes a cell, such as proteins or
other elements one can measure experimentally. The maturity is supposed to be a continuous variable
ranging from $m=0$ to $m=1$ in the two phases.

Cells enter the proliferating phase with age $a=0$. We suppose that they can be lost with a rate
$\gamma$. Every proliferating cell with maturity $m$ divides, at mitosis, with a rate
$\kappa=\kappa(m,a)$. We assume that the mapping $m\mapsto \kappa(m,a)$ is continuous, the mapping
$a\mapsto\kappa(m,a)$ is piecewise continuous, and that there exist $\underline{\tau}\geq0$ such
that
\begin{displaymath}
\kappa(m,a)=0 \qquad \textrm{ if } a\leq\underline{\tau}, \ m\in[0,1],
\end{displaymath}
and $0<\overline{\tau}<\infty$ such that
\begin{displaymath}
\int_0^{\overline{\tau}}\kappa(m,a)=+\infty \qquad \textrm{ for } m\in[0,1].
\end{displaymath}
This latter assumption describes the fact that a proliferating cell which has not died must divide
before it reaches the age $a=\overline{\tau}$. Therefore, the age of a proliferating cell ranges
from $a=0$ to $a=\overline{\tau}$.

Notice that $\underline{\tau}$, which represents the minimal time a cell has to spend in the
proliferating phase before it divides, can vanish. When $\underline{\tau}>0$, cells cannot divide
immediately after their entrance in the proliferating phase. This is probably the most realistic
situation since fast proliferations are usually observed in pathological cases.

All cells are supposed to age with unitary velocity, that means
\begin{displaymath}
\frac{da}{dt}=1,
\end{displaymath}
and to mature with a velocity $V(m)$, which is assumed to be the same in the two phases and satisfies $V(0)=0$, $V(m) > 0$ for $m \in (0,1]$, $V\in C^{1}[0,1]$ and
\begin{displaymath}
\int_{0}^{m} \frac{d\theta}{V(\theta)} = + \infty, \quad \textrm{ for } m \in (0,1].
\end{displaymath}
Since $\int_{m_{1}}^{m_{2}} \frac{d\theta}{V(\theta)}$ is the time required for a cell with maturity $m_{1}$ to reach a maturity $m_{2}$, this last condition on $V$ implies that a cell with maturity zero is unable to mature. For example, this condition is satisfied if $V(m)\underset{m\to 0}{\sim}\alpha_{1}m^{p}$, with $\alpha_{1}>0$ and $p\geq 1$.

After the division, each daughter cell enters immediately the resting phase with age $a=0$. If the maturity of the mother cell at the cytokinesis point was $m$, then we suppose that the two daughter cells have maturity $g(m)$ at birth. It is natural to assume that g is an increasing and continuous function from $[0,1]$ into itself. We also suppose that $g\in C^{1}[0,1)$ and $g(m)<m$ for $m \in (0,1)$. Furthermore, we suppose that $g^{-1}(m)=1$ if $m>g(1)$.

A cell can stay its entire life in the resting phase, so the age variable ranges from $a=0$ to
$a=+\infty$. We suppose that the resting cells can die at a rate $\delta$, which corresponds to
cellular differentiation, and that the reintroduction rate $\beta$, from the resting phase into the
proliferating phase, depends on the maturity and the total density of cells in the resting phase.

We denote by $p(t,m,a)$ and $n(t,m,a)$ respectively the population density in the proliferating and
the resting phase, at time $t\geq 0$, with maturity $m$ and age $a$. We define the total densities
of proliferating and resting cells at time $t$ and a maturity level $m$, respectively, by
\begin{displaymath}
P(t,m)=\int_{0}^{\overline{\tau}}p(t,m,a)da,
\end{displaymath}
and
\begin{displaymath}
N(t,m)=\int_{0}^{+\infty}n(t,m,a)da.
\end{displaymath}

Each phase is described by a partial differential equation. Hereafter, $p=p(t,m,a)$ and $n=n(t,m,a)$. The two conservation equations are
\begin{equation}\label{p(t,m,a)}
\frac{\partial p}{\partial t}+\frac{\partial p}{\partial a}+\frac{\partial (V(m)p)}{\partial m} =
-\left[\gamma(m)+\kappa(m,a)\right]p,
\end{equation}
\begin{equation}\label{n(t,m,a)}
\frac{\partial n}{\partial t}+\frac{\partial n}{\partial a}+\frac{\partial (V(m)n)}{\partial m} =
-\left[\delta(m)+\beta(m,N(t,m))\right]n.
\end{equation}

The functions $\gamma$ and $\delta$ are supposed to be positive and continuous on the interval $[0,1]$.

For the moment, $\beta$ is supposed to be a continuous and positive function. In fact, as we are studying a cellular population, we could consider only non-negative solutions ; hence, to give a correct definition of $\beta$, we will suppose that the mapping $x\mapsto \beta(m,x)$ is constant for $x\leq 0$, with $\beta(m,x)=\beta(m,0)$. For example, Mackey \cite{mackey1978} considered the function $\beta$ defined by
\begin{equation} \label{betalipschitz}
\beta(m,x)= \left\{ \begin{array}{ll}
\displaystyle\frac{a(m)}{x^{n}+b(m)}, & \textrm{ if } x\geq0, \\
\displaystyle\frac{a(m)}{b(m)}, & \textrm{ if } x<0,
\end{array} \right.
\end{equation}
for all $m\in[0,1]$, where $a$ and $b$ are two positive and continuous functions on $[0,1]$ and $n\geq 1$.

System (\ref{p(t,m,a)})-(\ref{n(t,m,a)}) is completed by boundary conditions modelling the flux between the resting phase and the proliferating phase.

The first condition,
\begin{equation}\label{bcp(t,m,a)}
p(t,m,0)=\beta (m,N(t,m)) N(t,m),
\end{equation}
represents the efflux of cells leaving the resting phase to the proliferating phase.

The second condition,
\begin{equation}\label{bcn(t,m,a)}
n(t,m,0)=2 (g^{-1})^{\prime}(m)\int_{\underline{\tau}}^{\overline{\tau}}
\kappa(g^{-1}(m),a)p(t,g^{-1}(m),a)da,
\end{equation}
describes the fact that, firstly, cells with maturity $m$ entering the resting phase come from
proliferating mother cells with a maturity $g^{-1}(m)$ and, secondly, that age of a dividing mother
cell is between $\underline{\tau}$ and $\overline{\tau}$ and this cell divides with a rate
$\kappa$.

We also consider the maturity flow $\pi_{s} : [0,1] \to [0,1]$, defined for $s \leq 0$, satisfying
the ordinary differential equation
\begin{displaymath}
\left\{ \begin{array}{rcll}
\displaystyle\frac{d\pi_{s}(m)}{ds} & = & V(\pi_{s}(m)),& s \leq 0, \\
\pi_{0}(m)       & = & m, &
\end{array} \right.
\end{displaymath}
for all $m\in[0,1]$, which represents the evolution of cells maturity to reach a maturity $m$ at time $0$ from a time $s \leq 0$. $\pi_{s}$ satisfies $\pi_{s}(0)=0$ and $\pi_{s}(m) \in (0,1]$ for $s\leq 0$ and $m \in (0,1]$.

One can remark that, for all $s\leq 0$ and $m\in[0,1]$,
\begin{displaymath}
\pi_{s}(m)=h^{-1}(h(m)e^{s}),
\end{displaymath}
where the continuous function $h :[0,1] \to [0,1]$ is defined by
\begin{displaymath}
h(m)=\left\{ \begin{array}{ll}
exp\bigg(-\displaystyle\int_{m}^{1} \frac{d\theta}{V(\theta)} \bigg),&\textrm{for } m \in (0,1],\\
0,& \textrm{for } m=0.
\end{array} \right.
\end{displaymath}

\bigskip{}

Integrating the equations (\ref{p(t,m,a)}) and (\ref{n(t,m,a)}) over the age, using the boundary
conditions (\ref{bcp(t,m,a)}) and (\ref{bcn(t,m,a)}) and the method of characteristics, we obtain
that $N(t,m)$ and $P(t,m)$ satisfy, for $t\geq\overline{\tau}$ and $m\in[0,1]$,
\begin{equation}\label{N(t,m)}
\begin{array}{c}
\displaystyle \frac{\partial }{\partial t}N(t,m)+\frac{\partial }{\partial m} (V(m)N(t,m)) =
-(\delta (m)+\beta (m,N(t,m)))N(t,m)
+2(g^{-1})^{\prime}(m)\times \vspace{1ex}\\
\displaystyle\displaystyle\int_{\underline{\tau}}^{\overline{\tau}} k(g^{-1}(m),a)
\xi(a,g^{-1}(m))\beta(\Delta(a,m),N(t-a,\Delta(a,m)))N(t-a,\Delta(a,m))da,
\end{array}
\end{equation}
\begin{equation}\label{P(t,m)}
\begin{array}{c}
\displaystyle \frac{\partial }{\partial t} P(t,m)+\frac{\partial }{\partial m} (V(m)P(t,m)) = -\gamma (m)P(t,m)+\beta (m,N(t,m))N(t,m)\\
\displaystyle -\int_{\underline{\tau}}^{\overline{\tau}}k(m,a)\xi(a,m)
\beta(\pi_{-a}(m),N(t-a,\pi_{-a}(m)))N(t-a,\pi_{-a}(m))da,
\end{array}
\end{equation}
where
\begin{displaymath}
\xi (t,m)=\exp\bigg\{-\int_{0}^{t} \gamma (\pi _{-s}(m))+V'(\pi _{-s}(m))ds\bigg\},
\end{displaymath}
and
\begin{displaymath}
k(m,a)=\kappa(m,a)\exp\left(-\int_0^a\kappa(\pi_{s-a}(m),s)ds\right).
\end{displaymath}

The mapping $\Delta : [0,+\infty)\times[0,1] \to \mathbb{R}^{+}$ is defined by
\begin{displaymath}
\Delta(s,m) = \pi_{-s}(g^{-1}(m)).
\end{displaymath}
$\Delta$ is non-increasing with respect to $s$ and non-decreasing with respect to $m$. Furthermore,
for $m \in (0,1]$,
\begin{equation}\label{delta}
\Delta(s,m) < m \quad \textrm{ if and only if } \quad s > \int_{m}^{g^{-1}(m)} \frac{d\theta}{V(\theta)}.
\end{equation}

We can notice that the solutions of Equation (\ref{P(t,m)}) depend on the solutions of Equation (\ref{N(t,m)}) whereas the converse is not true. In fact, if we know the behaviour of the solutions of Equation (\ref{N(t,m)}), which describes the process of replication, then we can deduce the behaviour of the solutions of Equation (\ref{P(t,m)}). Hence, we will focus our study on Equation (\ref{N(t,m)}), and particularly on the stability of its trivial solution.

Let $C[0,1]$ be the space of continuous functions on the interval $[0,1]$, endowed with the
supremum norm. We give, by using, for example, the semigroup theory, an integrated version of
Equation (\ref{N(t,m)}). First, we recall the expression of the semigroup generated by the linear
operator
\begin{displaymath}
A:\mathcal{D}(A) \subset C[0,1] \to C[0,1],
\end{displaymath}
defined by
\begin{displaymath}
\mathcal{D}(A) = \Big\{u \in C[0,1] \ ; \ u \in C^{1}(0,1], \ \lim_{m \to 0} V(m)u'(m) =0\Big\},
\end{displaymath}
and
\begin{displaymath}
Au(m)=\left\{\begin{array}{ll}
    -(\delta (m) +V'(m))u(m)-V(m)u'(m), & \textrm{if} \ m \in (0,1], \\
    -(\delta (0) +V'(0))u(0), & \textrm{if} \ m =0.
    \end{array} \right.
\end{displaymath}

We can easily show (one can refer to \cite{webb1996} for a similar operator) that the operator $A$
is the infinitesimal generator of a $C_{0}$-semi-groupe $(S(t))_{t \geq 0}$ defined for $\psi \in
C[0,1]$, $m \in [0,1]$ and $t \geq 0$, by
\begin{displaymath}
(S(t)\psi)(m)=\psi(\pi_{-t}(m))K(t,m),
\end{displaymath}
where
\begin{displaymath}
K(t,m)=exp\bigg\{-\int_{0}^{t} \delta (\pi_{-s}(m)) +V'(\pi_{-s}(m)) ds \bigg\},
\end{displaymath}
for all $m \in [0,1]$ and $t \geq 0$.

Hence, an integrated formulation (variation of constant formula) of Equation (\ref{N(t,m)}) is given by
\begin{equation} \label{mildformulation}
N(t,m) = \varphi(\overline{\tau},\pi_{-(t-\overline{\tau})}(m))K(t-\overline{\tau},m) +G(N)(t,m)-J(N)(t,m),
\end{equation}
if $t\geq \overline{\tau}$, and
\begin{equation} \label{icmf}
N(t,m) = \varphi(t,m),\ \  \textrm{ if } t\in[0,\overline{\tau}],
\end{equation}
for an initial condition $\varphi\in C([0,\overline{\tau}]\times[0,1])$, the space of continuous
functions on $[0,\overline{\tau}]\times[0,1]$, endowed with the supremum norm $||.||$, where
\begin{eqnarray*}
G(N)(t,m) = 2\!\!\int_{\overline{\tau}}^{t}\!\bigg(\!\int_{\underline{\tau}}^{\overline{\tau}}\!\! \zeta (a,\pi_{-(t-s)}(m))N(s-a,\Delta(a,\pi_{-(t-s)}(m)))\times  \\
\beta (\Delta(a,\pi_{-(t-s)}(m)),N(s-a,\Delta(a,\pi_{-(t-s)}(m)))  da \bigg) K(t-s,m) ds,
\end{eqnarray*}
where
\begin{displaymath}
\zeta(s,m)=(g^{-1})^{\prime}(m)k(g^{-1}(m),s)\xi(s,g^{-1}(m)),
\end{displaymath}
and
\begin{displaymath}
J(N)(t,m)=\int_{\overline{\tau}}^{t}K(t-s,m)\beta(\pi_{-(t-s)}(m),N(s,\pi_{-(t-s)}(m)))N(s,\pi_{-(t-s)}(m))ds,
\end{displaymath}
for $t \geq \overline{\tau}$ and $m \in [0,1]$.

One can also obtain, using the same technique, an integrated form of Equation (\ref{P(t,m)}).

Under classical Lipschitz conditions on the nonlinear part $x\mapsto x\beta(m,x)$, we can easily
show, using an iterative method and the method of steps, that, for each initial data $\varphi\in
C([0,\overline{\tau}]\times[0,1])$, the problem (\ref{mildformulation})-(\ref{icmf}) has a unique
continuous solution $N^{\varphi}$, which is nonnegative.

\section{Local stability} \label{sectionlocalstability}

We are interested, in this section, in the local stability of the trivial solution $N\equiv0$ of Equation (\ref{N(t,m)}) or its integrated version (\ref{mildformulation}). First, we recall a useful definition.
\begin{definition}
Let $\Omega$ be a subset of $C([0,\overline{\tau}]\times[0,1])$. We say that a solution $\overline{N}$ of Problem (\ref{mildformulation})-(\ref{icmf}), with an initial condition $\overline{\varphi}$, is {\em locally stable on the set $\Omega$\/} if, for all $\varepsilon>0$, there exists $\mu=\mu(\varepsilon)>0$ such that, if $\varphi\in\Omega$ and $||\varphi-\overline{\varphi}||<\mu$, then, for all $(t,m)\in[\overline{\tau},+\infty)\times[0,1]$, $|N^{\varphi}(t,m)-\overline{N}(t,m)|<\varepsilon$.\\
This solution $\overline{N}$ is said to be {\em locally exponentially stable on the set $\Omega$\/}
if there exists $\varepsilon>0$ such that, if $\varphi\in\Omega$ and $||\varphi -\overline{\varphi}
||<\varepsilon$, then there exist $c>0$ and $d>0$ such that, for all
$(t,m)\in[\overline{\tau},+\infty)\times[0,1]$,
\begin{equation} \label{caracexpstab}
|N^{\varphi}(t,m)-\overline{N}(t,m)|\leq ce^{-d(t-\overline{\tau})}.
\end{equation}
We will say {\em locally stable\/} (respectively {\em locally exponentially stable\/}) if
$\Omega=C([0,\overline{\tau}]\times[0,1])$.
\end{definition}

In this section, we make the assumption that the nonlinear function $x\mapsto x\beta(m,x)$, defined
for $x\in\mathbb{R}$, is Lipschitz continuous in a neighborhood of zero, for all $m\in[0,1]$, that
is, there exist $L\geq 0$ and $\varepsilon >0$ such that,
\begin{equation} \label{hypbeta}
|x\beta(m,x)-y\beta(m,y)| \leq L|x-y|,
\end{equation}
for all $|x|\leq\varepsilon$, $|y|\leq\varepsilon$ and $m\in[0,1]$. For example, if $\beta$ is
defined by (\ref{betalipschitz}), then the mapping $x\in\mathbb{R}\mapsto x\beta(m,x)$ is Lipschitz
continuous on $\mathbb{R}$ and a Lipschitz constant is given by $L=\sup_{m\in[0,1]}(a(m)/b(m))$.

In the sequel, $I$ and $\tilde{\zeta}$ denote the quantities
\begin{displaymath}
I:=\inf_{m\in[0,1]} (\delta(m)+V'(m)) \quad \textrm{ and } \quad
\tilde{\zeta}:=\int_{\underline{\tau}}^{\overline{\tau}}  \hat{\zeta}(a) da,
\end{displaymath}
where $\hat{\zeta}(a):=\sup_{m\in[0,1]}|\zeta(m,a)|$ and $\zeta(m,a)$ is given in the previous
section. Then we can prove the following invariance result, which is a first step to prove our
local stability result.

\begin{proposition} \label{localstability}
Suppose that $L(2\tilde{\zeta}+1) < I$. If $\varphi\in C([0,\overline{\tau}]\times[0,1])$ satisfies
$||\varphi||\leq\varepsilon$, where $\varepsilon$ is given by (\ref{hypbeta}), then, the solution
$N^{\varphi}$ of Problem (\ref{mildformulation})-(\ref{icmf}) satisfies
\begin{displaymath}
|N^{\varphi}(t,m)|\leq \varepsilon,
\end{displaymath}
for all $(t,m)\in[\overline{\tau},+\infty)\times[0,1]$.
\end{proposition}
\begin{proof}
It is clear that the assumption in this proposition implies that $I>0$. Let $\varepsilon >0$ be given by (\ref{hypbeta}) and $\varphi\in C([0,\overline{\tau}]\times[0,1])$ such that $||\varphi||\leq\varepsilon$.\\
We define the sequence $(N_{n})_{n\in\mathbb{N}}$ by
\begin{equation} \label{N0def}
N_{0}(t,m)=\left\{ \begin{array}{ll}
\varphi(\overline{\tau},\pi_{-(t-\overline{\tau})}(m))K(t-\overline{\tau},m), &  \textrm{ if } t\geq\overline{\tau}, \\
\varphi(t,m), &\textrm{ if } t\in[0,\overline{\tau}],
\end{array} \right.
\end{equation}
and, for $n\geq 1$,
\setlength\arraycolsep{2pt}
\begin{equation} \label{Nndef}
N_{n}(t,m)=\left\{ \begin{array}{ll}
N_{0}(t,m)+G(N_{n-1})(t,m)-J(N_{n-1})(t,m), & \textrm{ if } t\geq\overline{\tau}, \\
\varphi(t,m), &\textrm{ if } t\in[0,\overline{\tau}].
\end{array} \right.
\end{equation}
Since $\pi_{-(t-\overline{\tau})}(m) < m \leq g(1)$, then $|\varphi(\overline{\tau},\pi_{-(t-\overline{\tau})}(m))| \leq \varepsilon$. Moreover,
\begin{displaymath}
|K(t-\overline{\tau},m)|\leq e^{-I(t-\overline{\tau})},
\end{displaymath}
therefore
\begin{equation} \label{hypN0}
|N_{0}(t,m)| \leq \varepsilon e^{-I(t-\overline{\tau})} \leq \varepsilon.
\end{equation}
Let $T>\overline{\tau}$ be given. We show by induction that for all $n\in\mathbb{N}$ and all
$(t,m)\in[\overline{\tau},T]\times[0,1]$,
\begin{equation} \label{hyprecurrence}
|N_{n}(t,m)| \leq \varepsilon.
\end{equation}
For $n=0$, we just showed it in (\ref{hypN0}). Suppose that (\ref{hyprecurrence}) is true for $n\geq 0$. Then
\begin{displaymath}
|N_{n+1}(t,m)| \leq |N_{0}(t,m)|+|G(N_{n})(t,m)|+|J(N_{n})(t,m)|,
\end{displaymath}
with
\begin{displaymath}
|G(N_{n})(t,m)|\leq 2\tilde{\zeta}\varepsilon L\int_{\overline{\tau}}^{t} e^{-I(t-s)}ds,
\end{displaymath}
and
\begin{displaymath}
|J(N_{n})(t,m)|\leq \varepsilon L\int_{\overline{\tau}}^{t}e^{-I(t-s)}ds.
\end{displaymath}
So
\begin{eqnarray*}
|N_{n+1}(t,m)| & \leq & \varepsilon \bigg( e^{-I(t-\overline{\tau})} +L(1+2\tilde{\zeta})\int_{\overline{\tau}}^{t}e^{-I(t-s)}ds \bigg) \\
        & \leq & \varepsilon \bigg( e^{-I(t-\overline{\tau})} +L(1+2\tilde{\zeta})\frac{1-e^{-I(t-\overline{\tau})}}{I} \bigg).
\end{eqnarray*}
Set
\begin{displaymath}
\alpha=\frac{L(2\tilde{\zeta}+1)}{I}.
\end{displaymath}
Noting that $\alpha < 1$, we obtain
\begin{eqnarray*}
|N_{n+1}(t,m)| & \leq  & \varepsilon \Big( e^{-I(t-\overline{\tau})}(1-\alpha) + \alpha \Big), \\
        & \leq  & \varepsilon \big( (1-\alpha) + \alpha \big).
\end{eqnarray*}
That is, finally,
\begin{displaymath}
|N_{n+1}(t,m)| \leq  \varepsilon.
\end{displaymath}
Then, (\ref{hyprecurrence}) is true for every $n\in\mathbb{N}$.\\
Remarking that the sequence $(N_{n})_{n\in\mathbb{N}}$ converges uniformly to $N^{\varphi}$ for
$(t,m)\in[\overline{\tau},T]\times[0,1]$ (this can be proved using the method of steps), we
conclude that
\begin{displaymath}
|N^{\varphi}(t,m)|\leq \varepsilon.
\end{displaymath}
for all $t\in[\overline{\tau},T]$ and $m\in[0,1]$. This result is true for every
$T>\overline{\tau}$. So, it is true for $t\geq\overline{\tau}$ and the proof of the proposition is
ended.
\end{proof}

The result of invariance proved in Proposition \ref{localstability} is fundamental to obtain the local exponential stability of the trivial solution of Problem (\ref{mildformulation}). This result is presented in the following theorem.
\begin{theorem} \label{localexpstability}
If we suppose that $L(2\tilde{\zeta}+1) < I$, then the trivial solution $N\equiv 0$ of Problem
(\ref{mildformulation}) is locally exponentially stable.
\end{theorem}
\begin{proof}
Let $\varepsilon>0$ be given by (\ref{hypbeta}) and $\varphi\in C([0,\overline{\tau}]\times[0,1])$ such that $||\varphi||\leq\varepsilon$, and let $T>\overline{\tau}$ be given.\\
We consider the sequence defined by (\ref{N0def})-(\ref{Nndef}). From the proof of Proposition
\ref{localstability}, it follows that for all $n\in\mathbb{N}$ and
$(t,m)\in[\overline{\tau},T]\times[0,1]$,
\begin{displaymath}
|N_{n}(t,m)|\leq \varepsilon.
\end{displaymath}
Since the mapping
\begin{displaymath}
\rho\in[0,I]\mapsto \frac{I-\rho}{1+2\tilde{\zeta}e^{\rho\overline{\tau}}}
\end{displaymath}
is continuous on $[0,I]$ and reaches its maximum for $\rho = 0$, there exists $\rho\in(0,I)$ such that
\begin{displaymath}
L < \frac{I-\rho}{1+2\tilde{\zeta}e^{\rho\overline{\tau}}} < \frac{I}{2\tilde{\zeta}+1}.
\end{displaymath}
Using the same arguments as in the proof of Proposition \ref{localstability}, we obtain
\begin{displaymath}
|N_{0}(t,m)| \leq \varepsilon e^{-I(t-\overline{\tau})} \leq \varepsilon e^{-\rho(t-\overline{\tau})} \leq\varepsilon.
\end{displaymath}
The estimate
\begin{displaymath}
|N_{1}(t,m)-N_{0}(t,m)| \leq |G(N_{0})(t,m)|+|J(N_{0})(t,m)|,
\end{displaymath}
with
\begin{eqnarray*}
|J(N_{0})(t,m)| & \leq & L\int_{\overline{\tau}}^{t}|N_{0}(s,\pi_{-(t-s)}(m))|e^{-I(t-s)} ds, \\
        & \leq & \varepsilon L\int_{\overline{\tau}}^{t}e^{-I(t-s)}e^{-\rho(s-\overline{\tau})}ds, \\
        & \leq & \varepsilon Le^{-It}e^{\rho\overline{\tau}}\int_{\overline{\tau}}^{t}e^{(I-\rho)s}ds,
\end{eqnarray*}
and
\begin{eqnarray*}
|G(N_{0})(t,m)| & \leq & 2L\int_{\overline{\tau}}^{t}\bigg(\int_{\underline{\tau}}^{\overline{\tau}} \hat{\zeta}(a) |N_{0}(s-a,\Delta(a,\pi_{-(t-s)}(m)))| da \bigg) e^{-I(t-s)}ds,\\
     & \leq & 2\varepsilon L\int_{\overline{\tau}}^{t} e^{-I(t-s)}\bigg(\int_{\underline{\tau}}^{\overline{\tau}} \hat{\zeta}(a) e^{-\rho(s-a-\overline{\tau})} da \bigg) ds, \\
        & \leq & 2\tilde{\zeta}\varepsilon L\int_{\overline{\tau}}^{t} e^{-I(t-s)}e^{-\rho(s-2\overline{\tau})}ds, \\
        & \leq & 2\tilde{\zeta}\varepsilon Le^{-It}e^{2\rho\overline{\tau}}\int_{\overline{\tau}}^{t}e^{(I-\rho)s}ds,
\end{eqnarray*}
yields to
\begin{displaymath}
|N_{1}(t,m)-N_{0}(t,m)| \leq  \varepsilon Le^{-It}e^{\rho\overline{\tau}}\Big[
1+2\tilde{\zeta}e^{\rho\overline{\tau}} \Big]\int_{\overline{\tau}}^{t}e^{(I-\rho)s}ds.
\end{displaymath}
Note that
\begin{eqnarray*}
e^{-It}e^{\rho\overline{\tau}}\int_{\overline{\tau}}^{t} e^{(I-\rho)s} ds & = & e^{-It}e^{\rho\overline{\tau}}\frac{e^{(I-\rho)t}-e^{(I-\rho)\overline{\tau}}}{I-\rho}, \\
    & = & e^{(\rho-I)t}e^{\rho(\overline{\tau}-t)}\frac{e^{(I-\rho)t}-e^{(I-\rho)\overline{\tau}}}{I-\rho}, \\
    & = & e^{-\rho(t-\overline{\tau})}\frac{1-e^{(I-\rho)(\overline{\tau}-t)}}{I-\rho}, \\
    & \leq & \frac{1}{I-\rho}e^{-\rho(t-\overline{\tau})}.
\end{eqnarray*}
Consequently
\begin{displaymath}
|N_{1}(t,m)-N_{0}(t,m)| \leq  \varepsilon L\theta e^{-\rho(t-\overline{\tau})},
\end{displaymath}
where
\begin{displaymath}
\theta=\frac{1+2\tilde{\zeta}e^{\rho\overline{\tau}}}{I-\rho}.
\end{displaymath}
We show by induction that, for all $n\geq 1$ and all $(t,m)\in[\overline{\tau},T]\times[0,1]$,
\begin{equation} \label{hyprecurrence2}
|N_{n}(t,m)-N_{n-1}(t,m)| \leq  \varepsilon (L\theta)^{n} e^{-\rho(t-\overline{\tau})}.
\end{equation}
Suppose that (\ref{hyprecurrence2}) is true for $n\geq1$ ; then, for all
$(t,m)\in[\overline{\tau},T]\times[0,1]$,
\begin{displaymath}
|N_{n+1}(t,m)-N_{n}(t,m)| \leq  |(G(N_{n})-G(N_{n-1}))(t,m)|+|(J(N_{n})-J(N_{n-1}))(t,m)|,
\end{displaymath}
with
\setlength\arraycolsep{2pt}
\begin{eqnarray*}
|(G(N_{n})-G(N_{n-1}))(t,m)| & \leq & 2L\int_{\overline{\tau}}^{t} \bigg(\int_{\underline{\tau}}^{\overline{\tau}}\hat{\zeta}(a) |N_{n}(s-a,\Delta(a,\pi_{-(t-s)}(m))) \\
            &   & -N_{n-1}(s-a,\Delta(a,\pi_{-(t-s)}(m)))| da \bigg) e^{-I(t-s)} ds, \\
            & \leq & 2\varepsilon L^{n+1}\theta^{n}\!\!\int_{\overline{\tau}}^{t}\!\! e^{-I(t-s)}\bigg(\int_{\underline{\tau}}^{\overline{\tau}}\!\! \hat{\zeta}(a) e^{-\rho(s-a-\overline{\tau})}da\bigg) ds, \\
            & \leq & 2\tilde{\zeta}\varepsilon L^{n+1}\theta^{n}\int_{\overline{\tau}}^{t} e^{-I(t-s)} e^{-\rho(s-2\overline{\tau})} ds, \\
            & \leq & 2\tilde{\zeta}\varepsilon L^{n+1}\theta^{n}e^{2\rho\overline{\tau}}e^{-It}\int_{\overline{\tau}}^{t} e^{(I-\rho)s}ds,
\end{eqnarray*}
and
\begin{eqnarray*}
|(J(N_{n})-J(N_{n-1}))(t,m)| & \leq & L\int_{\overline{\tau}}^{t}|N_{n}(s,\pi_{-(t-s)}(m)) \\
            &   & \ -N_{n-1}(s,\pi_{-(t-s)}(m))| e^{-I(t-s)} ds, \\
            & \leq & \varepsilon L^{n+1}\theta^{n}\int_{\overline{\tau}}^{t}e^{-I(t-s)} e^{-\rho(s-\overline{\tau})} ds, \\
            & \leq & \varepsilon L^{n+1}\theta^{n}e^{\rho\overline{\tau}}e^{-It}\int_{\overline{\tau}}^{t}e^{(I-\rho)s} ds.
\end{eqnarray*}
Hence, finally, we obtain
\begin{eqnarray*}
|N_{n+1}(t,m)-N_{n}(t,m)| & \leq & \varepsilon L^{n+1}\theta^{n}e^{\rho\overline{\tau}}e^{-It}\Big[ 1+2\tilde{\zeta}e^{\rho\overline{\tau}} \Big]\int_{\overline{\tau}}^{t}e^{(I-\rho)s}ds, \\
            & \leq & \varepsilon (L\theta)^{n+1} e^{-\rho(t-\overline{\tau})}.
\end{eqnarray*}
This means that (\ref{hyprecurrence2}) is true for $n\geq 1$.\\
Then, we have, for all $(t,m)\in[\overline{\tau},T]\times[0,1]$,
\begin{eqnarray*}
|N_{n}(t,m)| & \leq & \sum_{i=1}^{n} |N_{i}(t,m)-N_{i-1}(t,m)| + |N_{0}(t,m)|, \\
    & \leq & \sum_{i=1}^{n} \varepsilon (L\theta)^{i} e^{-\rho(t-\overline{\tau})} +  \varepsilon e^{-\rho(t-\overline{\tau})}, \\
    & \leq & \varepsilon e^{-\rho(t-\overline{\tau})}\Big(1+L\theta+ \cdots + (L\theta)^{n} \Big).
\end{eqnarray*}
Since $L\theta < 1$, we conclude that, for all $t\in[\overline{\tau},T]$ and $m\in[0,1]$,
\begin{displaymath}
|N^{\varphi}(t,m)|\leq \frac{\varepsilon}{1-L\theta}e^{-\rho(t-\overline{\tau})}.
\end{displaymath}
This result is true for any $T>\overline{\tau}$. So, it is true for $t\geq\overline{\tau}$, and the
proposition is proved.
\end{proof}

\section{Global stability of the population of primitive cells} \label{stabilitydde}

We are interested, in this section, in the behaviour of the population of immature cells. That is, the cell population with maturity $m=0$.

Let $N^{\varphi}$ be a solution of Problem (\ref{N(t,m)})-(\ref{icmf}) for an initial condition $\varphi\in\mathcal{D}(A)$. We know that
\begin{displaymath}
\lim_{m\to 0} V(m)\frac{\partial}{\partial m}N^{\varphi}(t,m) = 0.
\end{displaymath}
Then, the function $x$ defined, for all $ t\geq 0$, by
\begin{displaymath}
x(t)=N^{\varphi}(t,0),
\end{displaymath}
is a solution of the following delay differential equation
\begin{equation} \label{dde}
x'(t) = -(I_{0}+\beta_{0}(x(t)))x(t) + 2\int_{\underline{\tau}}^{\overline{\tau}} Z(a)\beta_{0}(x(t-a))x(t-a)da,
\end{equation}
for $t\geq \overline{\tau}$, and
\begin{equation} \label{icdde}
x(t)=\psi(t):=\varphi(t,0), \ \ \textrm{for } t\in[0,\overline{\tau}],
\end{equation}
where
\begin{displaymath}
\left\{ \begin{array}{l}
I_{0}=\delta(0)+V'(0),\\
\beta_{0}(x)=\beta(0,x), \quad \textrm{for all } x\in\mathbb{R},\\
Z(a)=\zeta(a,0)=e^{-\nu a}k(0,a), \quad \textrm{with } \nu=\gamma(0)+V'(0).
\end{array} \right.
\end{displaymath}
It is clear that $I_{0}\geq 0$.

Since $\varphi\in C([0,\overline{\tau}]\times[0,1])$, then $\psi\in C[0,\overline{\tau}]$, where
$C[0,\overline{\tau}]$ is the space of all continuous functions on $[0,\overline{\tau}]$. Under
classical Lipschitz conditions on the function $\beta_{0}$, one can prove the existence and
uniqueness of solutions of Equation (\ref{dde}) (see Hale and Lunel \cite{halelunel}).

As in \cite{mackey1978}, it is natural to add the assumption that the function $\beta_{0}$ is decreasing for $x\geq 0$ and satisfies
\begin{displaymath}
\lim_{x\to +\infty} \beta_{0}(x)=0.
\end{displaymath}
A first property of the solutions of Problem (\ref{dde})-(\ref{icdde}) is given in the following proposition.
\begin{proposition}
Let $\psi\in C[0,\overline{\tau}]$. If $\psi$ is non-negative on $[0,\overline{\tau}]$, then the solution of Equation (\ref{dde}) for the initial condition $\psi$, denoted by $x^{\psi}$, is non-negative on $[\overline{\tau}, +\infty)$.
\end{proposition}
\begin{proof}
Suppose that the assumption of the proposition is false. First, suppose that there exist $t_{0}\in[\overline{\tau},\overline{\tau}+\underline{\tau})$ and $\varepsilon>0$ such that $ t_{0}+\varepsilon\leq\overline{\tau}+\underline{\tau}$, $x^{\psi}(t_{0})=0$ and $x^{\psi}(t)<0$ for $t\in(t_{0},t_{0}+\varepsilon)$. Then, for all $t\in[t_{0},t_{0}+\varepsilon)$,
\begin{displaymath}
-(I_{0}+\beta_{0}(x^{\psi}(t)))x^{\psi}(t) \geq 0,
\end{displaymath}
and
\begin{displaymath}
\int_{\underline{\tau}}^{\overline{\tau}} Z(a)\beta_{0}(x^{\psi}(t-a))x^{\psi}(t-a)da=\int_{\underline{\tau}}^{\overline{\tau}} Z(a)\beta_{0}(\psi(t-a))\psi(t-a)da \geq 0,
\end{displaymath}
because, for all $t\in[\overline{\tau},\overline{\tau}+\underline{\tau}]$ and all $a\in[\underline{\tau},\overline{\tau}]$, $t-a\in[0,\overline{\tau}]$. Hence, for $t\in[t_{0},t_{0}+\varepsilon)$, $(x^{\psi})'(t)\geq 0$, and we obtain a contradiction. So $x^{\psi}\geq 0$ on $[\overline{\tau},\overline{\tau}+\underline{\tau}]$.\\
Using the same arguments and the method of steps, we obtain the result for any
$t\geq\overline{\tau}$.
\end{proof}

As we are studying a cellular population, it is more interesting to consider non-negative initial conditions. Then we introduce the set
\begin{displaymath}
C^{+}=\big\{ \psi\in C[0,\overline{\tau}] \ ; \ \psi(t)\geq 0, \ \textrm{ for } \ t\in[0,\overline{\tau}] \big\}.
\end{displaymath}
The global stability of the population of primitive cells is given by the following result.
\begin{theorem} \label{ddestability}
Assume that $I_{0} > (2z-1)\beta_{0}(0)$ where $z=\int_{\underline{\tau}}^{\overline{\tau}}Z(a)da$. Then the trivial solution $x\equiv 0$ of Equation (\ref{dde}) is globally stable on the set $C^{+}$, that is, for all $\psi\in C^{+}$,
\begin{displaymath}
\lim_{t\to +\infty} x^{\psi}(t)=0.
\end{displaymath}
\end{theorem}
\begin{proof}
We prove the global stability of the trivial solution of (\ref{dde}) by using a Liapunov functional.\\
We define the functions $\lambda$ and $\Lambda$ by $\lambda(x)=x\beta_{0}(x)$ and $\Lambda(x)=\displaystyle\int_{0}^{x}\lambda(s)ds$, for all $x\geq0$. Let $H:C^{+} \to \mathbb{R}$ be the mapping defined, for all $\psi\in C^{+}$, by
\begin{displaymath}
H(\psi)=\Lambda(\psi(\overline{\tau}))+\int_{\underline{\tau}}^{\overline{\tau}}Z(a) \bigg( \int_{\overline{\tau}-a}^{\overline{\tau}}\lambda^{2}(\psi(\theta))d\theta \bigg) da.
\end{displaymath}
We show that $H$ is a Liapunov functional relative to Equation (\ref{dde}). We have
\begin{eqnarray}
\overset{\bullet}{H}(\psi) & = & \lim_{t \to 0} \frac{1}{t}(H(x_{t}^{\psi})-H(\psi)), \nonumber \\
    & = & \overset{\bullet}{\psi}(\overline{\tau})\lambda(\psi(\overline{\tau}))+\int_{\underline{\tau}}^{\overline{\tau}}Z(a) (\lambda^{2}(\psi(\overline{\tau}))-\lambda^{2}(\psi(\overline{\tau}-a)))da.\nonumber
\end{eqnarray}
Since
\begin{displaymath}
\overset{\bullet}{\psi}(\overline{\tau}) = -(I_{0}+\beta_{0} (\psi(\overline{\tau})))\psi(\overline{\tau}) + 2\int_{\underline{\tau}}^{\overline{\tau}} Z(a)\lambda(\psi(\overline{\tau}-a))da,
\end{displaymath}
then
\begin{eqnarray}
\overset{\bullet}{H}(\psi) & = & -[I_{0}+\beta_{0} (\psi(\overline{\tau}))]\beta_{0} (\psi(\overline{\tau}))\psi^{2}(\overline{\tau}) \nonumber \\
        &   & +\int_{\underline{\tau}}^{\overline{\tau}}Z(a) \big[\lambda^{2}(\psi(\overline{\tau}))+2\lambda(\psi(\overline{\tau}))\lambda(\psi(\overline{\tau}-a))-\lambda^{2}(\psi(\overline{\tau}-a))\big]da, \nonumber \\
        & = & -[I_{0}+\beta_{0} (\psi(\overline{\tau}))]\beta_{0} (\psi(\overline{\tau}))\psi^{2}(\overline{\tau}) + 2\int_{\underline{\tau}}^{\overline{\tau}}Z(a) \lambda^{2}(\psi(\overline{\tau}))da \nonumber \\
        &   & -\int_{\underline{\tau}}^{\overline{\tau}}Z(a) \big[\lambda(\psi(\overline{\tau}))-\lambda(\psi(\overline{\tau}-a))\big]^{2}da.  \nonumber
\end{eqnarray}
Hence,
\begin{displaymath}
\overset{\bullet}{H}(\psi) \leq  -[I_{0}-(2z-1)\beta_{0} (\psi(\overline{\tau}))]\beta_{0} (\psi(\overline{\tau}))\psi^{2}(\overline{\tau}).
\end{displaymath}
Let $u$ be the function defined, for $x \geq 0$, by
\begin{displaymath}
u(x)=r(x)\beta_{0}(x)x^{2},
\end{displaymath}
where, for all $x\geq 0$,
\begin{displaymath}
r(x)=I_{0}-(2z-1)\beta_{0}(x).
\end{displaymath}
$r$ satisfies $r(0)>0$ and $\lim_{x \to \infty} r(x)= I_{0}\geq 0$. Since $\beta_{0}$ is decreasing on $\mathbb{R}^{+}$, $r$ is a monotonous function. Therefore, $r$ is a positive function on $\mathbb{R}^{+}$.\\
Hence $u$ is non-negative on $[0,+\infty)$ and $u(x)=0$ if and only if $x=0$. Consequently, every
solution of Problem (\ref{dde})-(\ref{icdde}) with $\psi\in C^{+}$ tends to zero as $t$ tends to
$+\infty$.
\end{proof}

Combining this result with the result of local stability shown in Section \ref{sectionlocalstability}, we can obtain global stability for a class of solutions of Equation (\ref{N(t,m)}) or its integrated version (\ref{mildformulation}).

\section{Global stability}\label{globalstability}

In this section, we give a general criterion for global stability of Equation (\ref{mildformulation}). We first recall a definition.

\begin{definition}
Let $\Omega$ be a subset of $C([0,\overline{\tau}]\times[0,1])$. A solution $\overline{N}$ of
Problem (\ref{mildformulation})-(\ref{icmf}), with initial condition $\overline{\varphi}$, is said
to be {\em globally stable on the set $\Omega$\/} if, for all $\varphi\in\Omega$ and all
$m\in[0,1]$,
\begin{displaymath}
\lim_{t\to +\infty} |N^{\varphi}(t,m)-\overline{N}(t,m)|=0.
\end{displaymath}
This solution $\overline{N}$ is said to be {\em globally exponentially stable on the set
$\Omega$\/} if for all $\varphi\in\Omega$ there exist $c>0$ and $d>0$ such that, for all
$(t,m)\in[\overline{\tau},+\infty)\times[0,1]$,
\begin{displaymath}
|N^{\varphi}(t,m)-\overline{N}(t,m)|\leq ce^{-d(t-\overline{\tau})}.
\end{displaymath}
We will say  {\em globally stable\/} (respectively {\em globally exponentially stable\/}) if
$\Omega=C([0,\overline{\tau}]\times[0,1])$.
\end{definition}

From the result of local stability obtained in Section \ref{sectionlocalstability}, we are able to establish a first result about global stability for the trivial solution of Equation (\ref{mildformulation}). This result is presented in the following proposition.
\begin{proposition}\label{stabgloblip}
Suppose that the mapping $x\mapsto x\beta(m,x)$ is Lipschitz continuous on $[0,+\infty)$ for all
$m\in[0,1]$, with a Lipschitz constant $L$, such that
\begin{equation}\label{condstab}
L(2\tilde{\zeta}+1)<I.
\end{equation}
Then the trivial solution $N\equiv 0$ of Problem (\ref{mildformulation}) is globally exponentially stable.
\end{proposition}
\begin{proof}
We first show, as in the proof of Proposition \ref{localstability}, that for all $\varphi\in
C([0,\overline{\tau}]\times[0,1])$ and $(t,m)\in[0,+\infty)\times[0,1]$,
\begin{displaymath}
|N^{\varphi}(t,m)| \leq ||\varphi||.
\end{displaymath}
Then, we have (using the same arguments as in the proof of Theorem \ref{localexpstability}) that,
\begin{displaymath}
\lim_{t\to\infty}N^{\varphi}(t,m)=0 \quad \textrm{ exponentially},
\end{displaymath}
for all $\varphi\in C([0,\overline{\tau}]\times[0,1])$ and $m\in[0,1]$.
\end{proof}

The inequality
\begin{displaymath}
L(2\tilde{\zeta}+1)<I,
\end{displaymath}
is satisfied if $I$ or $\nu$ are large or if $L$ is small. Recall that $I$ and $\nu$ represent the
mortality rates in the resting and the proliferating phase and that $L$ is a bound of the quantity
of reintroduced cells. Then, the result of Proposition \ref{stabgloblip} represents a normal
behaviour of the cellular population.
\begin{corollary}
If the function $x\in\mathbb{R}\mapsto\beta(m,x)$ is differentiable and decreasing on $(0,+\infty)$
for all $m\in[0,1]$, then, under the assumption
\begin{displaymath}
(2\tilde{\zeta}+1) \sup_{m\in[0,1]}\beta(m,0) <I,
\end{displaymath}
the trivial solution $N\equiv 0$ of Problem (\ref{mildformulation}) is globally exponentially stable.
\end{corollary}
\begin{proof}
If the function $x\mapsto\beta(m,x)$ is differentiable and decreasing on $(0,+\infty)$ for all
$m\in[0,1]$, we obtain that
\begin{displaymath}
\frac{\partial}{\partial x}(x\beta(m,x)) \leq \beta(m,x) \leq \beta(m,0),
\end{displaymath}
for $x\in(0,+\infty)$ and $m\in[0,1]$. Then, the function $x\mapsto x\beta(m,x)$ is Lipschitz
continuous on $\mathbb{R}$ and the proposition \ref{stabgloblip} concludes the proof.
\end{proof}

For example, if $\beta$ is defined by (\ref{betalipschitz}) and if, in addition,
\begin{displaymath}
(2\tilde{\zeta}+1)\sup_{m\in[0,1]}\bigg(\frac{a(m)}{b(m)}\bigg) < I,
\end{displaymath}
then the trivial solution $N\equiv 0$ of Problem (\ref{mildformulation}) is globally exponentially stable.

Our next objective is to prove the global stability of solutions under more general conditions, less restrictive than the condition (\ref{condstab}) given in Proposition \ref{stabgloblip}.

We need the following assumption : there exists $C\geq0$ such that
\begin{equation} \label{integralefinie}
\int_{m}^{g^{-1}(m)} \frac{d\theta}{V(\theta)}\leq C,
\end{equation}
for all $m\in(0,1]$.

The biological interpretation of this condition is easy. That is, the time $\int_{m}^{g^{-1}(m)} \frac{d\theta}{V(\theta)}$ required for a daughter cell to reach the maturity $g^{-1}(m)$ of its mother is bounded. For example, if $V(m)\underset{m\to 0}{\sim}\alpha_{1}m$, with $\alpha_{1}>0$, then (\ref{integralefinie}) is equivalent to $g'(0)>0$.

Then we can prove the following fundamental theorem.

\begin{theorem}\label{smallmaturity}
Suppose that the condition (\ref{integralefinie}) is satisfied and that the function $x\mapsto
x\beta(m,x)$ is Lipschitz continuous for all $m\in[0,1]$. Set
\begin{displaymath}
\tau_{0} =\sup_{m>0} \Big( \int_{m}^{g^{-1}(m)} \frac{ds}{V(s)} \Big).
\end{displaymath}
Suppose that $N_{1}(t,m)$ and $N_{2}(t,m)$ are two solutions of Problem (\ref{mildformulation})-(\ref{icmf}) for the initial conditions $\varphi_{1}$ and $\varphi_{2}$ respectively. If
\begin{displaymath}
\underline{\tau} > \tau_{0},
\end{displaymath}
and if there exists $b \in (0,1]$ such that
\begin{displaymath}
\varphi_{1}(t,m)=\varphi_{2}(t,m),
\end{displaymath}
for $m \in [0,b]$ and $t \in [0,\overline{\tau}]$, then there exists $\overline{t} > \overline{\tau}$ such that
\begin{displaymath}
N_{1}(t,m)=N_{2}(t,m),
\end{displaymath}
for $m \in [0,1]$ and $t \geq \overline{t}$.
\end{theorem}

\begin{proof}
Since $\underline{\tau} > \tau_{0}$, then, from (\ref{delta}), $\Delta(a,m)<m$ for all
$a\in[\underline{\tau},\overline{\tau}]$.

We first show, using the Gronwall's lemma and a method of steps, that $N_{1}(t,m)=N_{2}(t,m)$ for
$m\in[0,b]$ and $t\geq 0$.

Let $ m \in [0,b]$ and $t\geq\overline{\tau}$ (if $t\in[0, \overline{\tau}]$, the result is
obvious). Then
\begin{displaymath}
\pi_{-(t-\overline{\tau})}(m) < m \leq b.
\end{displaymath}
So
\begin{displaymath}
\varphi_{1}(\overline{\tau},\pi_{-(t-\overline{\tau})}(m))=\varphi_{2}(\overline{\tau},\pi_{-(t-\overline{\tau})}(m)).
\end{displaymath}
Let $t \in [\overline{\tau},\overline{\tau}+\underline{\tau}]$. Then,
\begin{displaymath}
\Delta(a,\pi_{-(t-s)}(m)) < \pi_{-(t-s)}(m) < m \leq b,
\end{displaymath}
and consequently
\begin{displaymath}
G(N_{1})(t,m)=G(N_{2})(t,m).
\end{displaymath}
As the mapping $x\mapsto x\beta(m,x)$ is Lipschitz continuous and $K$ is bounded, we can write
\begin{eqnarray}
|N_{1}(t,m)-N_{2}(t,m)| & = & |J(N_{1})(t,m)-J(N_{2})(t,m)| \nonumber \\
            & \leq & c \int_{\overline{\tau}}^{t} |N_{1}(s,\pi_{-(t-s)}(m))-N_{2}(s,\pi_{-(t-s)}(m))| ds, \nonumber
\end{eqnarray}
where $c$ is a nonnegative constant. Then, using the Gronwall's lemma, we obtain that, for all
$t\in[\overline{\tau},\overline{\tau}+\underline{\tau}]$ and $m\in[0,b]$,
\begin{displaymath}
N_{1}(t,m)=N_{2}(t,m).
\end{displaymath}
Using the method of steps, we obtain this result for $t\geq\overline{\tau}$.

We define the mapping $\alpha:[0,1]\to[0,1]$ by
\begin{displaymath}
\alpha(m)=\Delta(\underline{\tau},m).
\end{displaymath}
Notice that $\alpha$ is an increasing function such that $\alpha(m) < m$ for all $m \in (0,1]$ and
$\alpha(g(1))=h^{-1}(e^{-\underline{\tau}})$. From (\ref{delta}) and since
$\underline{\tau}>\tau_0$, $h^{-1}(e^{-\underline{\tau}})<g(1)$. Hence, $\alpha$ is invertible and
we denote by $\Lambda$ its inverse function on $[0,h^{-1}(e^{-\underline{\tau}})]$. We set
$\Lambda(m)=g(1)$ for $m \in [h^{-1}(e^{-\underline{\tau}}),g(1)]$. Since $\alpha$ is non-negative
and increasing on the interval $[0,1]$, $\Lambda$ is increasing on
$[0,h^{-1}(e^{-\underline{\tau}})]$.

We first suppose that $b\in(0,h^{-1}(e^{-\underline{\tau}})]$.

We consider the increasing sequences $(b_{n})_{n \in \mathbb{N}}$ and $(t_{n})_{n \in \mathbb{N}}$
defined by
\begin{displaymath}
\left\{ \begin{array}{rcll}
b_{n+1} &=& \Lambda(b_{n}),& \textrm{ for } n \in \mathbb{N},\\
b_{0} &=& b,&
\end{array} \right.
\end{displaymath}
and
\begin{displaymath}
\left\{ \begin{array}{rcll}
t_{n+1} &=& t_{n}+\overline{\tau}+\displaystyle\ln\Big(\frac{h(b_{n+1})}{h(b_{n})}\Big),& \textrm{ for } n \in \mathbb{N},\\
t_{0} &=& 0.&
\end{array} \right.
\end{displaymath}
We show by induction that, if $m\in[0,b_{n}]$ and $t\geq t_{n}$, then $N_{1}(t,m)=N_{2}(t,m)$.\\
Suppose the equality is satisfied for $n\geq 0$. Let $m \in [0,b_{n+1}]$ and $t\geq t_{n+1}$. Then $t\geq t_{n}+\overline{\tau}$. The solutions  $N_{i}$, $i=1,2$, of Problem (\ref{mildformulation}) can be rewritten as
\setlength\arraycolsep{2pt}
\begin{displaymath}
\begin{array}{rcl}
N_{i}(t,m) & = & N_{i}(t_{n}+\overline{\tau},\pi_{-(t-t_{n}-\overline{\tau})}(m)) K(t-t_{n}-\overline{\tau},m)-J(N_{i})(t,m) \\
       & + & 2\displaystyle\int_{t_{n}+\overline{\tau}}^{t}\Big(\int_{\underline{\tau}}^{\overline{\tau}} \zeta(\pi_{-(t-s)}(m),a)N_{i}(s-a,\Delta(a,\pi_{-(t-s)}(m))) \times   \\
       &   & \beta(\Delta(a,\pi_{-(t-s)}(m)),N_{i}(s-a,\Delta(a,\pi_{-(t-s)}(m))))da \Big) K(t-s,m) ds,  \\
\end{array}
\end{displaymath}
for $t \geq t_{n}+\overline{\tau}$. Since
\begin{displaymath}
t-t_{n}-\overline{\tau} \geq \ln\bigg(\frac{h(b_{n+1})}{h(b_{n})}\bigg),
\end{displaymath}
we have
\begin{displaymath}
e^{-(t-t_{n}-\overline{\tau})} \leq \frac{h(b_{n})}{h(b_{n+1})}.
\end{displaymath}
Hence
\begin{displaymath}
h(m)e^{-(t-t_{n}-\overline{\tau})} \leq h(m)\frac{h(b_{n})}{h(b_{n+1})}\leq h(b_{n+1})\frac{h(b_{n})}{h(b_{n+1})}=h(b_{n}),
\end{displaymath}
because $h$ is an increasing function. Consequently,
\begin{displaymath}
\pi_{-(t-t_{n}-\overline{\tau})}(m)=h^{-1}(h(m)e^{-(t-t_{n}-\overline{\tau})}) \leq b_{n}.
\end{displaymath}
Using the induction hypothesis, we obtain that
\begin{displaymath}
N_{1}(t_{n}+\overline{\tau},\pi_{-(t-t_{n}-\overline{\tau})}(m))=N_{2}(t_{n}+\overline{\tau},\pi_{-(t-t_{n}-\overline{\tau})}(m)).
\end{displaymath}
Finally, since $s-a \geq t_{n}$ and
\begin{displaymath}
\Delta(a,\pi_{-(t-s)}(m)) \leq \Delta(a,m) \leq \Delta(\underline{\tau},m) \leq \Delta(\underline{\tau},b_{n+1})= b_{n},
\end{displaymath}
we obtain that
\begin{displaymath}
|N_{1}(t,m)-N_{2}(t,m)| = |J(N_{1})(t,m)-J(N_{2})(t,m)|.
\end{displaymath}
Using the Gronwall's lemma, we obtain the equality for $n+1$.

Noting that there exists $M \in \mathbb{N}^{*}$ such that $b_{M} < g(1)=b_{M+1}$, we obtain, for
all $m \in [0,b_{M+1}]=[0,1]$ and $t\geq t_{M+1}$, that $N_{1}(t,m)=N_{2}(t,m)$.

We now suppose that $b\in(h^{-1}(e^{-\underline{\tau}}),g(1)]$. Then the above result is still
valid with $M=0$. Thus,
\begin{equation}\label{property2}
N_1(t,m)=N_2(t,m) \qquad \textrm{ for } m\in[0,g(1)] \textrm{ and } t\geq t_{M+1}, \
M\in\mathbb{N}.
\end{equation}

We finally suppose that $b\in(g(1),1]$.

Let $m\in[g(1),1]$ and $t\geq t_{M+1}+\\overline{\tau}$ be given. For
$s\in[t_{M+1}+\overline{\tau},t]$ and $a\in[\underline{\tau},\overline{\tau}]$, we have $s-a\geq
t_{M+1}$. Thus, if $\chi(-(t-s),m)\leq g(1)$,
\begin{displaymath}
\Delta(a,\chi(-(t-s),m))\leq \Delta(a,g(1))\leq h^{-1}(e^{-\underline{\tau}})<g(1)
\end{displaymath}
and, from (\ref{property2}),
\begin{displaymath}
N_{1}(s-a,\Delta(a,\chi(-(t-s),m)))=N_{2}(s-a,\Delta(a,\chi(-(t-s),m))),
\end{displaymath}
and if  $\chi(-(t-s),m)>g(1)$ then, from the definition of $\zeta$,
\begin{displaymath}
\zeta(\chi(-(t-s),m),a)=0.
\end{displaymath}
Moreover, since $h(m)\leq 1$ for $m\in[0,1]$, we obtain, for $m\in[g(1),1]$ and $t\geq
t_{M+1}+\overline{\tau}-\ln(h(g(1)))$,
\begin{displaymath}
\chi(-(t-t_{M+1}-\overline{\tau}),m)\leq h^{-1}(h(m)h(g(1))) \leq g(1).
\end{displaymath}
We deduce that
\begin{displaymath}
N_{1}\big(t_{M+1}+\overline{\tau},\chi(-(t-t_{M+1}-\overline{\tau}),m)\big)=N_{2}\big(t_{M+1}+\overline{\tau},\chi(-(t-t_{M+1}-\overline{\tau}),m)\big).
\end{displaymath}
Using once again the Gronwall's Lemma, we get
\begin{displaymath}
N_{1}(t,m)=N_2(t,m) \qquad \textrm{ for } m\in[g(1),1] \textrm{ and } t\geq
t_{M+1}+\overline{\tau}-\ln(h(g(1))), \ M\in\mathbb{N}.
\end{displaymath}
We set $\overline{t}=t_{M+1}+\overline{\tau}-\ln(h(g(1)))$ and the theorem is proved.
\end{proof}

In the sequel, we suppose that Condition (\ref{integralefinie}) is satisfied, that the function
$x\mapsto x\beta(m,x)$ is Lipschitz continuous for all $m\in[0,1]$ and that $\underline{\tau} >
\tau_{0}$.

\begin{corollary}\label{extinction}
Let $\varphi \in C([0,\overline{\tau}]\times[0,1])$. If there exists $b \in (0,1]$ such that
\begin{displaymath}
\varphi(t,m)=0,
\end{displaymath}
for $m \in [0,b]$ and $t \in [0,\overline{\tau}]$, then there exists $\overline{t} > \overline{\tau}$ such that
\begin{displaymath}
N^{\varphi}(t,m)=0,
\end{displaymath}
for $m\in[0,1]$ and $t \geq \overline{t}$.
\end{corollary}

The result in Corollary \ref{extinction} describes the influence of stem cells (which are small maturity cells) on the total population. If the time of replication is large enough ($\underline{\tau} > \tau_{0}$), then a destruction of the population of stem cells affects all the system and the total population of cells dies out after a finite time $\overline{t}$. Biologically, this case corresponds to the aplastic anemia, a disease which destroys or injures stem cells, preventing the division.

Let $\varphi \in C([0,\overline{\tau}]\times[0,1])$ and $b\in(0,1]$. We define the function
$\varphi_{b} : [0,\overline{\tau}]\times[0,1] \to \mathbb{R}$ by
\begin{equation}
\varphi_{b}(t,m)= \left\{ \begin{array}{ll}
\varphi(t,m), & \textrm{ if } m\in[0,b], \\
\varphi(t,b), & \textrm{ if } m\in[b,g(1)].
\end{array} \right.
\end{equation}
It is clear that $\varphi_{b} \in C([0,\overline{\tau}]\times[0,1])$ for all $b\in(0,1]$.

To prove our main result on the global stability, we need two preliminary lemmas.

\begin{lemma}\label{Mackeyrudnicki1}
Let $\varphi \in C([0,\overline{\tau}]\times[0,1])$. Then, for all $T>\overline{\tau}$ and
$\varepsilon>0$, there exists $b\in(0,1]$ such that
\begin{displaymath}
|N^{\varphi_{b}}(t,m)-N^{\varphi}(t,0)|<\varepsilon,
\end{displaymath}
for $t\in[T,T+\overline{\tau}]$ and $m\in[0,1]$.
\end{lemma}

\begin{proof}
Let $\varphi \in C([0,\overline{\tau}]\times[0,1])$, $T>\overline{\tau}$ and $\varepsilon>0$ be
given. Using the continuity of the solutions of Equation (\ref{mildformulation}), we deduce that
there exists $\delta=\delta(\varepsilon)>0$ such that, if the function $\psi\in
C([0,\overline{\tau}]\times[0,1])$ satisfies
\begin{displaymath}
|\psi(t,m)-\varphi(t,0)|<\delta,
\end{displaymath}
for $(t,m)\in[0,\overline{\tau}]\times[0,1]$, then
\begin{displaymath}
|N^{\psi}(t,m)-N^{\varphi}(t,0)|<\varepsilon,
\end{displaymath}
for $(t,m)\in[T,T+\overline{\tau}]\times[0,1]$. Moreover, $\varphi$ is continuous, so there exists
$b\in(0,1]$ such that
\begin{displaymath}
|\varphi(t,m)-\varphi(t,0)|<\delta,
\end{displaymath}
for all $(t,m)\in[0,\overline{\tau}]\times[0,b]$. Hence, we obtain that
\begin{displaymath}
|\varphi_{b}(t,m)-\varphi(t,0)|<\delta,
\end{displaymath}
for $(t,m)\in[0,\overline{\tau}]\times[0,1]$, and so,
\begin{displaymath}
|N^{\varphi_{b}}(t,m)-N^{\varphi}(t,0)|<\varepsilon,
\end{displaymath}
for $t\in[T,T+\overline{\tau}]$ and $m\in[0,1]$.
\end{proof}

\begin{lemma}\label{Mackeyrudnicki1bis}
Let $\varphi \in C([0,\overline{\tau}]\times[0,1])$. Then, for all $b\in(0,1]$, there exists
$\overline{t}\geq 0$ such that
\begin{displaymath}
N^{\varphi_{b}}(t,m)=N^{\varphi}(t,m),
\end{displaymath}
for $t\geq\overline{t}$ and $m\in[0,1]$.
\end{lemma}

\begin{proof}
From the definition of $\varphi_{b}$ and Theorem \ref{smallmaturity}, the lemma is obvious.
\end{proof}

Now, we are able to prove that the global stability of the solutions of (\ref{mildformulation}) depends on their local stability and the global stability of the population of immature cells.

\begin{theorem} \label{Mackeyrudnicki2}
Suppose that the trivial solution $N\equiv 0$ of Equation (\ref{mildformulation}) is locally stable (resp. locally exponentially stable). Then this trivial solution is globally stable (resp. globally exponentially stable) on the set
\begin{displaymath}
\Omega=\Big\{ \varphi\in C([0,\overline{\tau}]\times[0,1]) \ ; \ \lim_{t\to +\infty}
N^{\varphi}(t,0)=0 \Big\}.
\end{displaymath}
\end{theorem}

\begin{proof}
$N\equiv 0$ is a locally stable solution of Equation (\ref{mildformulation}). Then, for all
$\varepsilon>0$, there exists $\delta>0$ such that, if $\varphi\in
C([0,\overline{\tau}]\times[0,1])$ satisfies
\begin{displaymath}
||\varphi||<\delta,
\end{displaymath}
then
\begin{displaymath}
|N^{\varphi}(t,m)|<\varepsilon,
\end{displaymath}
for $t\geq\overline{\tau}$ and $m\in[0,1]$.\\
Let $\varepsilon>0$ be fixed and $\varphi \in\Omega$. Then $\lim_{t\to\infty} N^{\varphi}(t,0)=0$. So there exists $T>\overline{\tau}$ such that,
\begin{displaymath}
|N^{\varphi}(t,0)|<\delta/2,
\end{displaymath}
for $t\geq T$. Furthermore, from Lemma \ref{Mackeyrudnicki1}, there exists $b\in(0,1]$ such that
\begin{displaymath}
|N^{\varphi_{b}}(t,m)-N^{\varphi}(t,0)|<\delta/2,
\end{displaymath}
for $(t,m)\in[T,T+\overline{\tau}]\times[0,1]$. Hence,
\begin{displaymath}
|N^{\varphi_{b}}(t,m)|<\delta,
\end{displaymath}
for $(t,m)\in[T,T+\overline{\tau}]\times[0,1]$. Since Equation (\ref{mildformulation}) is
autonomous, $N^{\varphi_{b}}$ becomes an initial condition for (\ref{mildformulation}) on the
interval $[T,T+\overline{\tau}]$. So, we obtain that
\begin{displaymath}
|N^{\varphi_{b}}(t,m)|<\varepsilon,
\end{displaymath}
for $t\geq T+\overline{\tau}$ and $m\in[0,1]$.\\
From Lemma \ref{Mackeyrudnicki1bis}, there exists $\overline{t}\geq0$ such that
\begin{displaymath}
N^{\varphi_{b}}(t,m)=N^{\varphi}(t,m),
\end{displaymath}
for $t\geq\overline{t}$ and $m\in[0,1]$, so
\begin{displaymath}
|N^{\varphi}(t,m)|<\varepsilon,
\end{displaymath}
for $t\geq max\{T+\overline{\tau},\overline{t}\}$ and $m\in[0,1]$. Consequently,
\begin{displaymath}
\lim_{t\to\infty} N^{\varphi}(t,m)=0,
\end{displaymath}
for all $\varphi\in\Omega$ and $m\in[0,1]$.\\
If $N\equiv 0$ is locally exponentially stable, the conclusion is obtained in the same way.
\end{proof}

The theorem \ref{Mackeyrudnicki2} stresses the influence of immature cells in the behaviour of the total population. That is, if the population of immature cells is stable, it implies that the total population is also stable.

Remark that the assumption
\begin{displaymath}
L(2\tilde{\zeta}+1) < I,
\end{displaymath}
given in Proposition \ref{stabgloblip}, implies, from Theorem \ref{localexpstability}, the local
exponential stability of the trivial solution of Equation (\ref{mildformulation}). Furthermore,
since $\beta_{0}(0)\leq L$, $z\leq\tilde{\zeta}$ and $I\leq I_{0}$, we have
\begin{displaymath}
(2z-1)\beta_{0}(0) < I_{0}.
\end{displaymath}
Hence, it follows, from  Theorem \ref{ddestability}, that
\begin{displaymath}
C^{+}\subset\Omega.
\end{displaymath}
Then, Proposition \ref{stabgloblip} becomes a simple consequence of Theorem \ref{Mackeyrudnicki2}.

We conclude with a general result of global stability, concerning the entire system of blood production.

If the trivial solution $N\equiv 0$ of Equation (\ref{mildformulation}) is globally stable on a set
$\Omega$, then the trivial solution $(N,P)\equiv(0,0)$ of the integrated version of Problem
(\ref{N(t,m)})-(\ref{P(t,m)}) is globally stable on the set $\Omega\times
C([0,\overline{\tau}]\times[0,1])$.

\end{document}